\documentclass[conference]{IEEEtran}

\IEEEoverridecommandlockouts
\usepackage{cite}

   \usepackage[pdftex]{graphicx}
   \graphicspath{{./img/}}
   \DeclareGraphicsExtensions{.pdf,.jpeg,.png}

\usepackage{amsfonts}
\usepackage{amssymb}
\usepackage[cmex10]{amsmath}
\usepackage[utf8]{inputenc}
\usepackage{textcomp}
\usepackage[hidelinks]{hyperref}
\usepackage{multirow}
\usepackage{booktabs,colortbl}
\usepackage{algorithmicx}
\usepackage{algorithm}
\usepackage{algcompatible}
\usepackage{multirow}
\usepackage{threeparttable}
\usepackage{pgfplotstable}

\begin{document}

\title{A Low-Rank Rounding Heuristic for Semidefinite Relaxation of Hydro Unit Commitment Problems}

\author{%

\IEEEauthorblockN{M.~Paredes and L.~S.~A.~Martins}
\IEEEauthorblockA{IBM Research\\
São Paulo, Brazil\\
\{mparedes,leonardo.martins\}@br.ibm.com}
\and
\IEEEauthorblockN{S.~Soares }
\IEEEauthorblockA{School of Electrical and Computer Engineering\\
University of Campinas, Brazil\\
dino@cose.fee.unicamp.br}

\thanks{This work was supported by the Brazilian National Council for Science and Technology (CNPq) and São Paulo Research Foundation (FAPESP).}%

}
% \IEEEoverridecommandlockouts
%\IEEEpubid{\makebox[\columnwidth]{978-1-5386-1379-5/17/\$31.00~
%\copyright2017
%IEEE \hfill} \hspace{\columnsep}\makebox[\columnwidth]{ }}

\maketitle
%ABSTRACT
\begin{abstract}
Hydro unit commitment is the problem of maximizing water use efficiency while minimizing start-up costs in the daily operation of multiple hydro plants, subject to constraints on short-term reservoir operation, and long-term goals.
A low-rank rounding heuristic is presented for the semidefinite relaxation of the mixed-integer quadratic-constrained formulation of this problem.
In addition to limits on reservoir and generator operation, transmission constraints are represented by an approximate AC power flow model.
In our proposed method, the mathematical program is equivalently formulated as a QCQP problem solved by convex relaxation based on semidefinite programming, followed by a MILP solution of undefined unit commitment schedules.
Finally, a rank reduction procedure is applied.
Effectiveness of the proposed heuristic is compared to branch-and-bound solutions for numerical case studies of varying sizes of the generation and transmission systems.
\end{abstract}

%INDEX TERMS

%INTRODUCTION
\section{Introduction}
The economic operation of hydro-dominated power systems is implicitly associated with efficient use of water resources for power generation.
The nonlinear characteristics of hydropower production models and transmission constraints pose computational challenges.
Most often such challenges are overcome by simplification of either or both models, as evidenced in the literature.
For instance, in~\cite{Sifuentes2007488} a formulation with AC power flow equations is presented at the expense of a linear representation of hydro production functions, whereas DC power flow formulations are more commonly proposed for the solution of large-scale instances of short-term hydro-thermal coordination problems, e.g.~\cite{291004}.
Because such linear representations of the production function invariably ignore head variation effects on productivity, linear piece-wise formulations have been presented over the years~\cite{1137622}.
Furthermore, hydro-electric unit commitment adds complexity to the problem due to the combinatorial nature of unit configuration optimization.
Lagrangian relaxation and dynamic programming~\cite{918302} were once common strategies in the literature for problem resolution.

Although it has been shown~\cite{IIASA} that load demand can be formulated as a non-anticipative variable, in this paper a deterministic model is considered.
Moreover, water inflows could also be modeled as stochastic variables, although such formulation is dependent on the availability of reliable data in the considered time resolution.

In the literature, short-term operation scheduling of hydrothermal power systems is commonly decomposed  into two properly coordinated subproblems.
In~\cite{589675}, the hydro unit commitment (HUC) problem is solved by means of a network flow model using priority-list-based dynamic programming with reservoir aggregation.
The strategy in~\cite{4562139} is based on mixed-integer linear programming (MILP) are proposed as an alternative to the so-called ``curse of dimensionality'' assailing dynamic programming methods.

A quadratic formulation with hydro unit efficiency functions is presented in~\cite{1626389} for individual representation of hydropower stations in a Lagrangian relaxation framework with a packing technique.
The HUC problem is then transformed into a linear master problem with construction of a nonlinear integer subproblem and consideration of forbidden operation zones.
This particular decomposition helps to reduce the duality gap in every iteration if compared to classical Lagrangian relaxation.
Subproblems are solved by a Quasi-Newton algorithm, implementing sequential quadratic programming (SQP). Lagrangian-based heuristics~\cite{1178813} can also approximate the solution of the HUC problem by linearizing the hydropower model with inner minimization of thermal subproblems.

In~\cite{Catalao2010904} a nonlinear HUC model is presented, where the objective is a nonlinear function of water release with linear constraints.
Integer variables are used to represent model characteristics related to forbidden zones and unit startup and shutdown.
Quadratic programming (QP) and barrier methods are used to solve the problem.

In this paper a quadratic model~\cite{6919349} of the hydroelectric production function is used, such that water discharge can be formulated as a function of active power generation, i.e:
\begin{equation} \label{eq:waterd}
 \hat{q}_{t,h}(.) = \sum_{u=1}^{N_{U_{t,h}}} x_{t,h,u} \left( \alpha_{h,u} p_{t,h,u}^{2} + \beta_{h,u} p_{t,h,u} + \gamma_{h,u} \right),
\end{equation}
where $r$ represents water discharge, and $\alpha,\beta,\gamma$ are coefficients of the water discharge curve for given active power output $p$ and unit configuration $u$ at hydro plant $h$ and hour $t$.
The choice of unit configurations is determined by discrete variables:
\begin{equation}
  x_{t,h,u} \in \{0,1\}, \qquad \forall (t,h,u)
\end{equation}
subject to hourly configuration uniqueness constraints:
\begin{equation}
  \sum_{u = 1}^{N_{U_{t,h}}} x_{t,h,u} = 1
\end{equation}
and unit startup availability between successive hours :
\begin{equation}
  \sum_{u=1}^{N_{U_{t,h}}} \left(u \cdot x_{t,h,u} \right) - \sum_{u=1}^{N_{U_{t-1,h}}} \left(u \cdot x_{t-1,h,u} \right) \leq y_{t,h}
\end{equation}
where $N_{U_{t,h}}$ is the number of unit configurations available for commitment, and $y_{t,h} \geq 0$ is the number of units scheduled for startup.
In order to solve this mixed-integer, quadratically-constrained quadratic program (MIQCQP) in a computationally efficient manner, we initially solve a semidefinite programming-based (SDP) convex relaxation~\cite{anjos:2012} of the problem with AC power flow formulation, as originally proposed in~\cite{6919349}.
This convex SDP relaxation has been previously shown to guarantee valid lower bounds with relatively small gaps.
In this paper we propose a rounding heuristic that, after the initial relaxation, solves a MILP subproblem to determine the optimal unit configuration of the remaining undefined discrete variables for the active power generation dispatches obtained from the solution to the SDP relaxation.
Once the active power dispatches and unit commitments are obtained, a re-optimization of the AC power flow using SDP is performed, followed by a rank-reduction procedure for improved solution.

%SOLUTION METHODOLOGY
\section{Solution Methodology}

\subsection{Problem formulation}
\label{sec:solution:qcqp}

It is based on the quadratic relaxation of discrete variables of the day-ahead HUC problem with hourly discretization~\cite{6919349}, such that:
\begin{equation}
\label{eq:solution:qcqp}
\begin{array}{crclcl}
\displaystyle \min_{\mathbf{x}} & \mathbf{x}^{\intercal}\mathbf{C}\mathbf{x}\\
\textrm{subject to} & \mathbf{x}^{\intercal}\mathbf{A}_{i}\mathbf{x} &=& a_{i} & & i=1,\ldots,m\\
    & \mathbf{x}^{\intercal} \mathbf{B}_{j}\mathbf{x} & \leq & b_{j}  & & j=1,\ldots,n\\
     & \mathbf{x} & \succeq & \mathbf{0}\\
\end{array}
\end{equation}
where $\mathbf{C}$ represents water use and unit startup costs, $\mathbf{A}_{i}$ represents active and reactive power balance, long-term reservoir operation goals, unit configuration uniqueness constraints, and fixed voltage at the slack bus.
Limits on active and reactive power generation, reservoir storage, as well as power flow and bus voltages are represented in $\mathbf{B}_{j}$.
Problem variables are represented by
\begin{displaymath}
  \mathbf{x} = \left( 1,x, y, p, q, e, f \right)^{\intercal}
\end{displaymath}
where, in addition to unit configuration $x$ and startup control variables $y$, active and reactive power, as well as real and imaginary parts of bus voltages in Cartesian coordinates, represented by $p$, $q$, $e$, and $f$, respectively.

The next sections describe the solution methodology  proposed in this paper as depicted in Algorithm~\ref{alg:solution}.

\begin{algorithm}[t]
  \begin{center}
    \begin{algorithmic}[1]
      \STATE Solve SDP (P1) to obtain $\tilde{p}, \tilde{x}$
      \STATE $\Omega_{1}(\tilde{\mathbf{x}}) \leftarrow \{(t,h,u) : \tilde{x}_{t,h,u} \notin \{ 0,1\} \}$
      \STATE $\Omega_{2}(\tilde{\mathbf{x}}) \leftarrow \{(t,h,u) : \tilde{x}_{t,h,u} = 1 \}$
      \STATE Let $\tilde{p}_{t,h} \leftarrow \sum_{u=1}^{N_{U_{t,h}}} \tilde{p}_{t,h,u} \cdot \tilde{x}_{t,h,u}$
      \STATE Solve MILP (P2) for $\tilde{p},\tilde{x}\in\Omega_{2}(\tilde{\mathbf{x}})$ to obtain $x^{*} \in \{0,1\}$
      \STATE Solve the OPF \eqref{eq:solution:acpf} for $ x^{*}$ to obtain $\hat{p}, \hat{q}, \hat{e}, \hat{f}$
      \STATE Apply rank reduction procedure to obtain $p^{*}, q^{*}, e^{*}, f^{*}$
      \STATE \textbf{Return} $p^{*},x^{*},q^{*},y^{*},e^{*},f^{*}$
    \end{algorithmic}
  \end{center}
  \caption{Solution methodology}
  \label{alg:solution}
\end{algorithm}

%SEMIDEFINITE RELAXATION
\subsection{Semidefinite relaxation}
\label{sec:solution:sdp}
The SDP relaxation~\cite{nesterov:2000} of the QCQP~\eqref{eq:solution:qcqp} model is given by:
\begin{eqnarray*}
\begin{array}{lrclrcl}
  \displaystyle (P1) \min_{\mathbf{X},  \mathbf{V}, \mathbf{y}} & \multicolumn{3}{c}{ \mathbf{C} \bullet \mathbf{X}  + \mathbf{c}^{\intercal} \mathbf{y}} \\
  \textrm{s. t.} & \mathbf{p} &=& \mathcal{P}(\mathbf{X})  + \mathcal{F}_{\mathbf{p}}(\mathbf{V}) ;& \mathbf{1} &=& \mathcal{A}(\mathbf{X}) \\
    & \mathbf{q} &=& \mathcal{Q}(\mathbf{X})+ \mathcal{F}_{\mathbf{q}}(\mathbf{V}) ;&  \mathbf{g} &=& \mathcal{G}(\mathbf{X})  \\
    & \mathbf{0} &\geq& \overline{\mathcal{P}}(\mathbf{X}) & \mathbf{0} &\geq& \overline{\mathcal{Q}}(\mathbf{X})  \\
    & \mathbf{r} &\geq& \mathcal{R}(\mathbf{X})  & \mathbf{s} &=& \mathcal{S}(\mathbf{V})\\
    & \overline{\mathbf{f}} &\geq& \mathcal{F}(\mathbf{V}) &  \overline{\mathbf{v}} &\geq& \mathcal{J}(\mathbf{V})  \\
    & \mathbf{y} &\geq& \mathcal{E}(\mathbf{X}) & 0 & \preceq & \mathbf{X}, \mathbf{V}\\
    & 1 & = & \textrm{rank}(\mathbf{X}) = \textrm{rank}(\mathbf{V})
    & \mathbf{0} &\leq& \mathbf{y}

\end{array}
\label{eq:SDPrelax}
\end{eqnarray*}
where the rank constraints on $\mathbf{X}$ and $\mathbf{V}$ are relaxed, and the Frobenius product is represented by the symbol $\bullet$.
\begin{equation*}
\mathbf{x}_{t,h,u}= \left[\Delta P_{t,h,u},x_{t,h,u},\Delta Q_{t,h,u} \right]^{\intercal} ; \; \mathbf{y}=  \left\{ y_{t,h,u } \right\}     \;  \forall t,h,u
\end{equation*}
\begin{equation*}
\mathbf{X} = \sum_{\left(t,h,u \right)  }  \boldsymbol{\xi}_{t,h,u }\boldsymbol{\xi}_{ t,h,u }^{\intercal} \otimes \mathbf{x}_{t,h,u}\mathbf{x}_{t,h,u}^{\intercal} ; \; \mathbf{V}=\sum_{t=1}^{T} \boldsymbol{\xi}_{t}\boldsymbol{\xi}_{t}^{\intercal} \otimes \mathbf{v}_{t}\mathbf{v}_{t}^{\intercal}
\end{equation*}
\begin{equation*}
\mathbf{v}_{t}=\left[ e_{t,1},\ldots, e_{t, N_{B}},f_{t,1},\ldots,f_{t,N_{B}} \right]^{\intercal} \forall t=1,\ldots,T
\end{equation*}
where $\boldsymbol{\xi}_{(t,h,u)}$ is a vector on the Euclidean space with appropriate dimensions (in this case $\mathbb{R}^{T \cdot N_{U}}$), with value 1 in the corresponding $(t,h,u)$ position.
Matrix $\mathbf{C}$ represents water cost coefficients:
\begin{align}
\mathbf{C}&= \sum_{t,h,u } \lambda_{h} \cdot  \mathbf{\widehat{C}}_{t,h,u} \\
\mathbf{\widehat{C}}_{t,h,u}&= \boldsymbol{\xi}_{t,h,u }\boldsymbol{\xi}_{t,h,u }^{\intercal} \otimes\left[ \begin{array}{ccc} \alpha_{h,u}  &  \frac{\widehat{\beta}_{h,u} }{2}  & 0\\ \frac{\widehat{\beta}_{h,u} }{2}  & \widehat{\gamma}_{h,u} & 0 \\ 0 & 0 & 0 \end{array} \right] \forall t,h,u
\end{align}
Startup costs are in vector $\mathbf{c} \in \mathbb{R}^{ \dim (\mathbf{y})}$, with elements
$c_{t,h,u}= \delta_{h} \; ; \; \forall t,h,u$.
Active power generation constraints are represented by linear mapping $\mathcal{P} \left( \cdot \right)$ , where, for each $\mathbf{P}_{i,t} \in \mathbb{R}^{ \dim(\mathbf{x}) \times \dim(\mathbf{x})}$ $\forall i= 1,\ldots,N_{B}$ and $t = 1,\ldots,T$, we have:
\begin{equation}
\mathbf{P}_{i,t} =  \sum_{h \in \Psi_{i}}  \sum_{u=1}^{Nu_{t,h}}\boldsymbol{\xi}_{t,h,u }\boldsymbol{\xi}_{t,h,u  }^{\intercal} \otimes \left[ \begin{array}{ccc} 0 &  \frac{1}{2} & 0  \\ \frac{1}{2}  & \underline{P}_{h,u} & 0 \\ 0 & 0 & 0 \end{array} \right]
\end{equation}
%\begin{equation}
%\mathbf{Q}_{i,t} =  \sum_{h \in \Psi_{i}}  \sum_{u=1}^{Nu_{t,h}}\boldsymbol{\xi}_{t,h,u }\boldsymbol{\xi}_{t,h,u }^{\intercal} \otimes \left[ \begin{array}{ccc} 0 &  0 & 0  \\ 0  & \underline{Q}_{h,u} & \frac{1}{2} \\ 0 & \frac{1}{2} & 0 \end{array} \right]
%\end{equation}
where $\underline{P}_{h,u}$ is the minimum active power generation, Linear mapping $\mathcal{Q} \left( \cdot \right)$ is similarly constructed.
In order to represent the active and reactive power injection on every bus we use $\mathcal{F}_{\mathbf{p}}\left( \cdot \right)$ and $\mathcal{F}_{\mathbf{q}} \left( \cdot \right)$, with elements $\mathbf{F_{p}}_{i,t}$ and $\mathbf{F_{q}}_{i,t} \in \mathbb{R}^{2\cdot N_{B} \times 2\cdot N_{B}}$ .
\begin{equation}
\mathbf{F_{p}}_{i,t} =-\sum_{j \in \Omega_{i}}\boldsymbol{\xi}_{t}\boldsymbol{\xi}_{t}^{\intercal} \otimes \mathbf{Y}_{i,j} \; ; \; \mathbf{F_{q}}_{i} =-\sum_{j \in \Omega_{i}} \boldsymbol{\xi}_{t}\boldsymbol{\xi}_{t}^{\intercal} \otimes\overline{\mathbf{Y}}_{i,j}
\end{equation}
\begin{equation}
\mathbf{Y}_{i,j} = \frac{1}{2} \left[ \begin{array}{cc}
G_{i,j}+G_{i,j}^{\intercal} & B_{i,j}^{\intercal}-B_{i,j} \\
 B_{i,j}-B_{i,j}^{\intercal} & G_{i,j}+G_{i,j}^{\intercal}
\end{array}
\right]
\end{equation}
\begin{equation}
\overline{\mathbf{Y}}_{i,j} = - \frac{1}{2} \left[ \begin{array}{cc}
B_{i,j}+B_{i,j}^{\intercal} & G_{i,j}^{\intercal}-G_{i,j} \\
 G_{i,j}-G_{i,j}^{\intercal} & B_{i,j}+B_{i,j}^{\intercal}
\end{array}
\right]
\end{equation}
\begin{equation*}
G_{i,j}=g_{i,j} \left( \boldsymbol{\xi}_{i}\boldsymbol{\xi}_{i}^{\intercal}- \boldsymbol{\xi}_{i}\boldsymbol{\xi}_{j}^{\intercal}  \right) \textrm{ and }  B_{i,j}=b_{i,j} \left( \boldsymbol{\xi}_{i}\boldsymbol{\xi}_{i}^{\intercal}- \boldsymbol{\xi}_{i}\boldsymbol{\xi}_{j}^{\intercal}  \right)
\end{equation*}
where $g_{i,j}$ and $b_{i,j}$ are conductance and susceptance, and vectors $\mathbf{p}$ and $\mathbf{q}$ contain active and reactive power loads.
The value of the module of the slack-bus voltage $V_{slack}$ is represented by $\mathcal{S} \left( \cdot \right)$, with matrices $\mathbf{S}_{t}  \in \mathbb{R}^{ 2 \cdot N_{B} \times 2 \cdot N_{B}} $,
\begin{equation}
\mathbf{S}_{t}=\boldsymbol{\xi}_{i_{slack}}\boldsymbol{\xi}_{i_{slack}}^{\intercal} \; ; \;\mathbf{S}_{T+t}=\boldsymbol{\xi}_{N_{B}+i_{slack}}\boldsymbol{\xi}_{N_{B}+i_{slack}}^{\intercal}
\end{equation}
and vector $ \mathbf{s} \in \mathbf{R}^{ 2 \cdot T}$, with elements $s_{t}=V_{slack} , s_{T+t}=0 , \; \forall t =1,\ldots,T$.
Active power flow is limited using linear mapping $\mathcal{F} \left( \cdot \right)$ where $\mathbf{F}_{l,t} \in \mathbb{R}^{ 2 \cdot N_{B} \times 2 \cdot N_{B}}$ is given by:
\begin{equation}
\mathbf{F}_{l,t}=-\mathbf{F}_{ N_{L}+ r,t}=\boldsymbol{\xi}_{t}\boldsymbol{\xi}_{t}^{\intercal} \otimes\mathbf{Y}_{i,j}.
\end{equation}
Vector $ \overline{\mathbf{f}} \in \mathbb{R}^{2 \cdot T \cdot N_{L}}$ of power flow limits has elements $\overline{f}_{l,t}=\overline{f}_{N_{L}+r,t}= \overline{F}_{i,j}$
where $\forall t=1,\ldots, T$, such that $(i,j) \in \mathcal{L}$.
Voltage variables at bus $i$ and time $t$ are limited in linear mapping $\mathcal{J}_{i,t}\left( \cdot \right)$ where matrices $\mathbf{J}_{i,t} \in \mathbb{R}^{ 2 \cdot N_{B}\cdot T \times 2 \cdot N_{B} \cdot T} $ are given by:
\begin{equation}
\mathbf{J}_{i,t} = -\mathbf{J}_{N_{B} +i,t} =\boldsymbol{\xi}_{t}\boldsymbol{\xi}_{t}^{\intercal} \otimes\left(\boldsymbol{\xi}_{i}\boldsymbol{\xi}_{i}^{\intercal} +\boldsymbol{\xi}_{N_{B}+i}\boldsymbol{\xi}_{N_{B}+i}^{\intercal} \right)
\end{equation}
and voltage limits ($\overline{V}_{i}$ and $\underline{V}_{i}$) are obtained from $\overline{\mathbf{v}} \in \mathbb{R}^{2 \cdot N_{B}}$
\begin{equation}
\overline{v}_{i,t}=\overline{V}_{i} \;\; ;\;\; \overline{v}_{N_{B} +i,t}=-\underline{V}_{i}
\end{equation}
$ \forall  t=1, \ldots , T \; ; \; i= \{ 1, \ldots , N_{B} \} \backslash i_{slack}$.

Active power generation target constraint at hydro plant $h$ is represented in linear mapping $\mathcal{G} \left( \cdot \right)$, where matrices $\mathbf{G}_{h} \in \mathbb{R}^{ \dim(\mathbf{x1}) \times \dim(\mathbf{x1})}$ are given by:
\begin{equation}
\mathbf{G}_{h} =   \sum_{t=1}^{T}  \sum_{u=1}^{Nu_{t,h}}\boldsymbol{\xi}_{t,h,u }\boldsymbol{\xi}_{t,h,u }^{\intercal} \otimes \left[ \begin{array}{ccc} 0 &  \frac{1}{2} & 0  \\ \frac{1}{2}  & \underline{P}_{h,u} & 0 \\ 0 & 0 & 0 \end{array} \right]
\end{equation}
$\forall h=\{ 1, \ldots, N_{H}  \} \backslash h_{slack}$,
where $\mathbf{g}$ contains active power generation targets obtained from longer-term hydropower scheduling models.
Active and reactive power generation limits are represented by linear mappings $\overline{\mathcal{P}} \left(\cdot \right)$  contain matrices $\overline{\mathbf{P}}_{t,h,u} \in \mathbb{R}^{ \dim(\mathbf{X}) \times \dim(\mathbf{X})}$:
\begin{equation}
\overline{\mathbf{P}}_{t,h,u } =  \boldsymbol{\xi}_{t,h,u }\boldsymbol{\xi}_{t,h,u}^{\intercal} \otimes \left[ \begin{array}{ccc} 0 &  \frac{1}{2} & 0  \\ \frac{1}{2}  & -\overline{ \Delta P}_{h,u} & 0 \\ 0 & 0 & 0 \end{array} \right]
\end{equation}
%\begin{equation}
%\overline{\mathbf{Q}}_{t,h,u} =  \boldsymbol{\xi}_{t,h,u}\boldsymbol{\xi}_{t,h,u}^{\intercal} \otimes \left[ \begin{array}{ccc} 0 &  0 & 0  \\ 0  & -%\overline{ \Delta Q}_{h,u} & \frac{1}{2} \\ 0 & \frac{1}{2} & 0 \end{array} \right]
%\end{equation}
This is similar for linear mapping $\overline{\mathcal{Q}} \left(\cdot \right)$.
Reservoir dynamics at hydro plant $h$ and time $t$ is given in $\mathcal{R}\left(\cdot \right)$ containing matrices $\mathbf{R}_{h,t} \in \mathbb{R}^{ \dim(\mathbf{X}) \times \dim(\mathbf{X}) }$, and $\mathbf{R}_{h,t}= -\mathbf{R}_{N_{H}+h,t} $ such that
\begin{equation*}
\mathbf{R}_{h,t}= \vartheta \sum_{i=1}^{t} \left( \sum_{u=1}^{Nu_{i,h}} \widehat{\mathbf{C}}_{i,h,u}  + \sum_{\overline{h} \in \Theta_{h} } \sum_{u=1}^{Nu_{i-\tau,\overline{h}}}  \widehat{\mathbf{C}}_{i-\tau,\overline{h},u}   \right)
\end{equation*}
and the vector  $\mathbf{r} \in \mathbb{R}^{2 \cdot T \cdot N_{H}}$ with elements:
\begin{equation}
r_{h,t}=\overline{\nu}_{h}  -\nu o_{h}  - \vartheta \sum_{i=1}^{t} \left( a_{h,i}-s_{h} + \sum_{\tilde{h} \in \Theta_{h}}  s_{\overline{h}}\right)
\end{equation}
such that $r_{ N_{H}+h,t} = \overline{\nu}_{h} -\underline{\nu}_{h} -r_{h,t}$, $\forall t=1,\ldots,T$ and $h=1,\ldots, N_{H}$, where $a_{h,i}$, $\nu o_{h}$, and $s_{h}$ represent given water inflow, initial reservoir volume and spillage, respectively.

Unit configuration status variables and their respective uniqueness constraints are represented in $\mathcal{A}\left(\cdot\right)$ where matrices $\mathbf{A}_{h,t}\in\mathbb{R}^{\dim \left(\mathbf{X} \right)}$ are given by:
\begin{equation}
\mathbf{A}_{h,t}= \sum_{u=1}^{Nu_{t,h}}  \boldsymbol{\xi}_{t,h,u} \boldsymbol{\xi}_{t,h,u}^{\intercal} \otimes    \left[\begin{matrix}
0 & 0 & 0 \\
0 & 1 & 0 \\
0 & 0 & 0
\end{matrix} \right]
\end{equation}
$\forall t=1,\ldots, T$ and $h=1,\ldots,N_{H}$.
Startup constraints uses $\mathcal{E} \left( \cdot \right)$ with $\mathbf{E}_{t,h} \in \mathbb{R}^{ \dim \left( \mathbf{x} \right) }$, such that:
\begin{equation*}
\mathbf{E}_{t,h} = \sum_{u=1}^{Nu_{t,h}} \left( \boldsymbol{\xi}_{t,h,u } \boldsymbol{\xi}_{t,h,u }^{\intercal} -  \boldsymbol{\xi}_{t-1,h,u } \boldsymbol{\xi}_{t-1,h,u}^{\intercal} \right)\otimes \left[
\begin{matrix}
0 & 0 & 0 \\
0 & u & 0 \\
0 & 0 & 0
\end{matrix}
\right]
\end{equation*}
$\forall t=1,\ldots, T$ and $h=1,\ldots,N_{H}$.

%SOLVE UNDEFINED UNIT COMMITMENT SCHEDULES
\subsection{Solving undefined unit commitment schedules}
A feasible active power dispatch $\tilde{p}$ is obtained from the solution to the SDP relaxation in (P1), which may contain undefined unit commitment schedules represented by the set:
\begin{equation}
  \Omega_{1}(\tilde{\mathbf{x}}) = \left\{ (t,h,u) : \tilde{x}_{t,h,u} \notin \left\{ 0,1\right\} \right\}
\end{equation}
A MILP problem is formulated at domain $x\in\Omega_{1}(\tilde{\mathbf{x}})$ with:
\begin{equation}
P_{t,h}^{*}= \sum_{u=1}^{Nu_{t,h}} \textbf{1}_{\left(t,h,u\right) \in \Omega_{1} \left(\mathbf{x} \right) } \cdot \left( \Delta P_{t,h}^{*} + \underline{P}_{h} \cdot x_{t,h}^{*} \right)
\end{equation}
such that:
\begin{equation*}
  \begin{array}{crclcl}
    \displaystyle (P2) \min_{x \in \Omega_{1}(\tilde{\mathbf{x}}), y} & \multicolumn{3}{c}{ \displaystyle \sum_{\substack{(t,h,u)  \\ \in \Omega_{1}(\tilde{\mathbf{x}})}} wc_{t,h} \resizebox{.09 \textwidth}{!}
{$(\Delta P_{t,h,u}^{*}, x_{t,h,u}^{*}) $} x_{t,h,u}  \displaystyle + \sum_{t=1,h=1}^{T,N_{H}} \delta_{h}  y_{t,h}}\\
    \textrm{subject to } &\displaystyle \sum_{u=1}^{Nu_{t,h}} \textbf{1}_{(t,h,u) \in \Omega_{1 }(\tilde{\mathbf{x}})}  x_{t,h,u}&=&1\\
\multicolumn{2}{l}{ \displaystyle \sum_{u=1}^{Nu_{t,h}} \left( \textbf{1}_{(t,h,u) \in \Omega_{1 }(\tilde{\mathbf{x}})}  u  x_{t,h,u} +  \textbf{1}_{(t,h,u) \in \Omega_{2 }(\tilde{\mathbf{x}})} u \right) }\\
  \multicolumn{2}{r}{  \displaystyle - \sum_{u=1}^{Nu_{t-1,h}} \left( \textbf{1}_{(t-1,h,u) \in \Omega_{1 }(\tilde{\mathbf{x}})} u  x_{t-1,h,u} \right. }  \\
 \multicolumn{2}{r}{ \displaystyle \left.   + \textbf{1}_{(t-1,h,u) \in \Omega_{2}(\tilde{\mathbf{x}})}  u  \right) } &\leq& y_{t,h}
  \end{array}
  \label{eq:solution:milp}
\end{equation*}
for $t = 1,\ldots,N_{T}$ and $h = 1,\ldots,N_{H}$, where $\textbf{1}_{w \in A}=\{1$, if $w \in A$, $0$ otherwise$\}$. Water use and unit start-up costs are represented by $wc(\cdot)$ and $ \delta_{h} $ respectively.
A feasible unit commitment schedule $x^{*}$ is obtained from solving (P2).

%AC POWER FLOW RE-OPTIMIZATION
\subsection{AC power flow re-optimization}
Since active power generation dispatches $\tilde{p}$ remain feasibly unchanged after solving (P2) for the undefined unit commitment schedules, it is necessary to guarantee feasibility of the updated water discharge values resulting from $x^{*}$ whilst minimizing costs with water use in a stripped down formulation of~\eqref{eq:SDPrelax}.
This is accomplished by solving the following SDP formulation of the optimal AC power flow problem to find $\hat{p},\hat{q},\hat{e},$ and $\hat{f}$:
\begin{equation}
  \begin{array}{crclcl}
    \displaystyle \min_{\mathbf{Z}} & \multicolumn{3}{l}{\widehat{\mathbf{C}} \bullet \mathbf{Z} } \\
    \textrm{subject to } & \widehat{\mathcal{A}}(\mathbf{Z})&=& \widehat{\mathbf{a}} & &  \\
      & \widehat{\mathcal{B}}(\mathbf{Z}) & \leq & \widehat{\mathbf{b}}  \\
      & \mathbf{Z} & \succeq & 0 \\
  \end{array}
  \label{eq:solution:acpf}
\end{equation}
where $\mathbf{Z}=\mathbf{z} \mathbf{z}^{T}$ is the outer product of $\mathbf{z}=(1,p,q,e,f)^{T}$, water use costs are represented by $\widehat{\mathbf{C}}$, $\widehat{\mathcal{A}}(\cdot)$ represents active and reactive power balance, long-term reservoir operation goals, and fixed voltage at the slack bus. Limits on active and reactive power generation, reservoir storage, as well as power flow and bus voltages are defined by $\widehat{\mathcal{B}}(\cdot)$.
Analogously to the initial relaxation, the rank-1 constraint on $\mathbf{Z}$ is relaxed.

%RANK REDUCTION PROCEDURE

\subsection{Rank reduction procedure}
\label{sec:rank}
In our proposed rounding heuristic we seek to further improve $\mathbf{Z}^{*}$ resulting from solving~\eqref{eq:solution:acpf} by applying rank reduction in an iterative procedure~\cite{1997geometry}, however, by updating matrix $\mathbf{Z}^{*}_{r}$ :
\begin{equation}
  \label{eq:rankup}%
  \mathbf{Z}^{*}_{r+1} \leftarrow \mathbf{Z}^{*}_{r} +\omega_{r}\cdot \mathbf{D}_{r}
\end{equation}
 Also for the current OPF SDP formulation the solution's rank is upper-bounded by the tree width of the power network plus one as demonstrated in~\cite{7065336}. The rank reduction procedure listed in Algorithm~\ref{alg:alg2} maintains feasibility as explained in the remainder of this section, and also converges up to a theoretical limit of the rank in the final solution~\cite{barvinok2002course}.
Moreover, because $\mathbf{Z}$ does not appear in the objective function, new solutions resulting from the rank reduction procedure maintain the value of the objective function. A central idea to this procedure is to calculate $\mathbf{D}_{r}$ such that at every iteration it decreases the rank of $\mathbf{Z}^{*}_{r}$ maintaining its feasibility, i.e.~$\mathbf{D}_{r}$ is be subject to:
\begin{algorithm}[t]
\begin{center}
\begin{algorithmic}[1]
\REQUIRE $\mathbf{Z}^{*}$,$ \{ \mathbf{Y}_{i} \forall i \in \{1,\ldots,N_{B} \}\}$, $\{ \mathbf{Y}_{i,j} \forall (i,j) \in \mathcal{L} \}$,  $\{\mathbf{E}_{l} \forall l \in \{1,\ldots,2N_{B}  \} \}$ and $\mathbf{Z}^{*}_{1}=\mathbf{Z}^{*}$
\FOR{$r\leftarrow 1,\ldots,\left( \text{rank}(\mathbf{Z}^{*})+\varsigma-1 \right) $}
\STATE Obtain Cholesky decomposition $\mathbf{Z}^{*}_{r}=\mathbf{R}_{r}\mathbf{R}_{r}^{\intercal}$
\STATE Obtain $\mathbf{\varepsilon}_{i}$ by get the basis of~\eqref{eq:feas2-Yi}-~\eqref{eq:feas2-El} .
\FOR{$i = 1,\ldots,\operatorname{nullity}(\mathbf{A})$ }
\STATE Calculate  $ \mathbf{S}_{r} = \operatorname{findSDP}(\mathbf{S}_{r},\operatorname{smat}(\mathbf{\varepsilon}_{i}))$
\ENDFOR
\STATE $\mathbf{D}_{r} \leftarrow - \mathbf{R}_{r}\mathbf{S}_{r}\mathbf{R}_{r}^{\intercal}$
\STATE $\mathbf{Z}^{*}_{r+1} \leftarrow \mathbf{Z}^{*}_{r} + \lambda_{\max}(\mathbf{S}_{r})  \cdot \mathbf{D}_{r}$
\IF{$\text{rank}(\mathbf{Z}^{*}_{r+1}) = 1$ \textbf{or} $\mathbf{S}_{r} = \mathbf{0} $}
\STATE \textbf{break}
\ENDIF
\ENDFOR
\end{algorithmic}
\end{center}
\caption{Rank reduction procedure.}\label{alg:alg2}
\end{algorithm}
\begin{align}
  \label{eq:feas-Yi}%
  \mathbf{D}_{r} \bullet \mathbf{Y}_{i} & = 0, \, i=1,\ldots,N_{B}  \\
  \label{eq:feas-Yij}%
  \mathbf{D}_{r} \bullet \mathbf{Y}_{i,j} & = 0, \, \forall\left(i,j\right) \in \mathcal{L} \\
  \label{eq:feas-El}%
  \mathbf{D}_{r} \bullet \mathbf{E}_{l} & = 0, \, l = 1,2,\ldots,2N_{B}
\end{align}
where $\mathbf{D}_{r} = \mathbf{R}_{r} \mathbf{S}_{r} \mathbf{R}_{r}^{\intercal}$  and matrix $\mathbf{R}_{r}$ are obtained from the Cholesky decomposition of $\mathbf{Z}^{*}_{r} = \mathbf{R}_{r}\mathbf{R}_{r}^{\intercal}$.
By the distributive property of the Frobenius product we obtain the following relationship:
\begin{equation}
\label{eq:distributive-property}%
 \mathbf{D}_{r}\bullet\mathbf{Y}_{i}  = \mathbf{R}_{r} \mathbf{S}_{r} \mathbf{R}_{r}^{\intercal} \bullet \mathbf{Y}_{i} = \mathbf{S}_{r} \bullet \mathbf{R}_{r}^{\intercal}  \mathbf{Y}_{i}  \mathbf{R}_{r} = 0.
\end{equation}
Similarly for the other constraints, we obtain a linear system:
\begin{align}
  \label{eq:feas2-Yi}%
  \operatorname{svec}\left(\mathbf{R}_{r}^{\intercal} \mathbf{Y}_{i}  \mathbf{R}_{r}\right)^{\intercal}\operatorname{svec}\left(\mathbf{S}_{r} \right)   & = 0, \, i=1,\ldots,N_{B}  \\
  \label{eq:feas2-Yij}%
 \operatorname{svec}\left( \mathbf{R}_{r}^{\intercal} \mathbf{Y}_{i,j} \mathbf{R}_{r} \right)^{\intercal}\operatorname{svec}\left( \mathbf{S}_{r}\right)& = 0, \, \left(i,j\right) \in \mathcal{L} \\
  \label{eq:feas2-El}%
\operatorname{svec}\left( \mathbf{R}_{r}^{\intercal} \mathbf{E}_{l}  \mathbf{R}_{r}\right)^{\intercal}\operatorname{svec}\left( \mathbf{S}_{r}\right)&= 0, \, l = 1,2,\ldots,2N_{B}
\end{align}
where $\operatorname{svec}(\cdot):\mathbb{R}^{n \times n} \rightarrow \mathbb{R}^{\frac{n(n+1)}{2}}$ is a linear transformation of a symmetric matrix into a vector.

We can represent the linear system in \eqref{eq:feas2-Yi}, \eqref{eq:feas2-Yij} and  \eqref{eq:feas2-El} as $\mathbf{A} \mathbf{s} = \mathbf{0}$, which is equivalent to say that $\mathbf{s}$ is in the null space of $\mathbf{A}$.
In order to solve this linear system SVD decomposition of matrix $\mathbf{A}=\mathbf{U \Sigma L}^{\intercal}$ is used.
If matrix $A$ is rank-deficient, then the null space is generated by the columns of $\mathbf{L}$ corresponding to the zero singular values on $\mathbf{\Sigma}$.
Then $\operatorname{svec}\left(\mathbf{S}_{r}\right)$ is generated by the basis of the null space of $\mathbf{A}$.
Finally, it is necessary to enforce positive semidefiniteness on $\mathbf{S}_{r}$ as the convex combination of the basis vectors of the null space of $\mathbf{A}$:
\begin{equation}\label{eq:feasibility}
\mathbf{S}_{r} = \sum_{i=1}^{\operatorname{nullity}(\mathbf{A}) } \operatorname{smat}(\mathbf{\varepsilon}_{i}) \sigma_{i}
\end{equation}
% then the following SDP problem must be solved:
%
%\begin{equation}\label{eq:rankred}
%\begin{array}{lrcll}
%  \displaystyle \min_{\mathbf{S}_{r},  \mathbf{z}} & \multicolumn{3}{c}{ z_{0}} \\
%    \textrm{s. t.}  &\displaystyle \mathbf{S}_{r} &=& \sum_{i=1}^{\operatorname{nullity}(\mathbf{A}) } \operatorname{smat}(\mathbf{\varepsilon}_{i}) z_{i} + \mathbf{I}z_{0} & \\
%    & 0 & \preceq & \mathbf{S}_{r}\\
%    & 0 &\leq& z_{0}\\
%\end{array}
%\end{equation}
where $\mathbf{\varepsilon}_{i}$ is a basis of the null space of $\mathbf{A}$.
The search for the semidefinite convex combination uses the algorithm in~\cite{crawford1983} to determine if is feasible to obtain a semidefinite matrix from a convex combination of two matrices.
This algorithm is therefore identified as $\operatorname{findSDP}(\cdot,\cdot) : \mathbb{R}^{n\times n} \times \mathbb{R}^{n\times n} \rightarrow \mathbb{R}^{n\times n}$.

% If problem~\eqref{eq:rankred} have no solution, then its not possible to reduce the rank of the matrix.
%In order to satisfy~\eqref{eq:distributive-property}, $\mathbf{S}_{r}$ is constructed as a $0$-$1$ matrix such that there is a one-valued element for every zero-valued element in the corresponding coordinate of matrices $\left( \mathbf{R}_{r}^{\intercal}  \mathbf{Y}_{i}  \mathbf{R}_{r} \right)$, $\left( \mathbf{R}_{r}^{\intercal}  \mathbf{Y}_{i,j}  \mathbf{R}_{r} \right)$, and $\left( \mathbf{R}_{r}^{\intercal}  \mathbf{E}_{l}  \mathbf{R}_{r} \right)$. This guarantees that~\eqref{eq:feas-Yi}, \eqref{eq:feas-Yij}, and~\eqref{eq:feas-El} are satisfied.

Then from~\eqref{eq:rankup} we obtain:
\begin{equation}
  \label{eq:rankred}%
  \mathbf{Z}^{*}_{r} + \omega_{r} \cdot \mathbf{D}_{r} = \mathbf{R}_{r} \left( \mathbf{I} +\omega_{r}\mathbf{S}_{r} \right) \mathbf{R}_{r}^{\intercal}.
\end{equation}
From~\eqref{eq:rankred} we conclude that the necessary condition for positive semidefiniteness of the new solution is given by:
\begin{equation}
  \left( \mathbf{Z}^{*}_{r} + \omega_{r}\cdot \mathbf{D}_{r} \right)   \succeq \mathbf{0} \Leftrightarrow \left( \mathbf{I} + \omega_{r}\cdot \mathbf{S}_{r} \right) \succeq \mathbf{0}.
\end{equation}
If $\mathbf{S}_{r}$ is indefinite, then the algorithm is unable to reduce its rank at the $r$-th iteration.
It should be noted, however, that in theory it should still be possible to reduce its rank in further iterations, given that $\mathbf{S}_{r}\neq\mathbf{S}_{r+1}$ because of the added perturbations, and also because the SDP formulation of the network constraints guarantees that at least one rank-1 solution exists.
In that sense, additional $\varsigma$ iterations are provisioned.
Parameter $\omega_{r}$, on the other hand, is a scale of perturbation matrix $\mathbf{D}_{r}$ that induces rank reduction of $\mathbf{Z}^{*}_{r}$ in at least one degree by reducing the rank of $\left( \mathbf{I}+ \omega_{r}\mathbf{S}_{r} \right)$.
This is achieved by obtaining the maximum eigenvalue of $\mathbf{S}_{r}$: $\omega_{r} = -\lambda_{\max}(\mathbf{S}_{r})$.

%NUMERICAL RESULTS
\section{Numerical Results}
The proposed solution methodology was implemented in MATLAB\textregistered~R2016.
The {SDPA}~\cite{sdpa} solver was used for solve SDP problems. The HUC MILP subproblem is solved using IBM\textregistered~CPLEX\textregistered~10.2.
A total of five numerical cases studies were analyzed for an increasing size of the test systems. Hydro plant data were obtained from the Brazilian power system operator for actual plants for a typical day.
Network data were obtained from modified versions of IEEE test systems.
Objective function results and gaps relative to the initial SDP relaxation are listed on Table~\ref{tab:results:costs} for the initial semidefinite relaxation (SDP), the B\&B algorithm solving a SDP relaxation in every node (this method is an exact method developed on \cite{6919349}), as well as the rounding heuristic (RH).
Table~\ref{tab:results:time} lists the computer time required by each of the procedures.
\begin{table}[t]
  \centering
  \begin{threeparttable}
    \caption{Objective Function and Relative Gap}
    \label{tab:results:costs}
    \setlength{\extrarowheight}{1pt}
    \begin{scriptsize}
      \begin{tabular}{lrrrrr}
        \hline
        \multirow{2}{*}{Case} & \multicolumn{3}{c}{Obj. Func. (\$)} & \multicolumn{2}{c}{Relative Gap (\$)}\\
        \cline{2-6}
         & \multicolumn{1}{c}{SDP} & \multicolumn{1}{c}{B\&B} & \multicolumn{1}{c}{RH} & \multicolumn{1}{c}{B\&B} & \multicolumn{1}{c}{RH}\\
        \hline
        3-GEN & 1,598,129 & 1,598,969 & 1,598,974 & 840 & 845\\
        IEEE-14 & 2,711,258 & 2,711,718 & 2,711,823 & 459  & 565\\
        IEEE-30 & 5,587,150 & 5,587,329 & 5,588,811 & 178 & 1,661\\
        IEEE-57 & 13,221,266 & 13,221,609 & 13,222,070 & 342 & 803\\
        IEEE-118 & 29,659,377 & ? & 29,703,639 & ? & 44,262\\
        \hline
      \end{tabular}
    \end{scriptsize}
  \end{threeparttable}
\end{table}

\begin{table}[t]
  \centering
  \begin{threeparttable}
    \caption{Computer Time Comparison}
    \label{tab:results:time}
    \setlength{\extrarowheight}{1pt}
    \begin{scriptsize}
      \begin{tabular}{lrrr}
        \hline
        \multirow{2}{*}{Case} & \multicolumn{3}{c}{Time (s)}\\
        \cline{2-4}
         & SDP & B\&B & RH\\
        \hline
        3-GEN & 0.78 & 45.64 & 22.13 \\
        IEEE-14 & 1.83  & 46.72 & 45.16 \\
        IEEE-30 & 10.00 & 560.594 & 20.32\\
        IEEE-57 & 33.03 & 10,475.46 & 64.55 \\
        IEEE-118 & 270.15 & $\infty$ & 1,009.97 \\
        \hline
      \end{tabular}
    \end{scriptsize}
  \end{threeparttable}
\end{table}
It can be observed that the solutions obtained by the B\&B and rounding heuristic algorithms were equivalent up to very small error margins. In the case of the IEEE 118-bus problem, the B\&B algorithm does not converge.
Results in  Table~\ref{tab:results:costs} suggest that the first relaxation is very close to the optimal solution obtained by B\&B , and consequently that it is possible to approximate the objective function with a small error margin, and, in the case of the proposed rounding heuristic, with reduced computing time. In all the cases the final unit commitment give us different result for B\&B and RH methods.
\section{Conclusion}
For the study cases presented the proposed rounding heuristic efficiently approximates the optimal solutions of different HUC problems. This is most likely due to the fact that the initial semidefinite relaxation already finds near-optimal solutions, as evidenced by the results obtained by the B\&B algorithm.
Furthermore, our rank reduction heuristic enhances the method presented in~\cite{1997geometry} by adding a procedure to identify the possibility to find a perturbation matrix that reduces in at least 1 the rank of the solution per iteration.
% TODO: Perhaps should remove the sentence below
% We should probably not leave the reader with open-ended comments in the conclusion - to highlight the computational efficiency of our procedure would require to elaborate what would be this SDP model
% Also, it is possible furmulate a feasibility SDP based on\eqref{eq:feasibility}, but numerical experiments performed suggests that solve this SDP is computationally expensive than use the search procedure.
It is possible to employ a more conservative rounding heuristic, as long as the tradeoff between computational burden and relative gap is considered.
However, we note that the observed performance of the proposed rounding heuristic is compatible with the results presented in~\cite{Goemans:1995:IAA:227683.227684}, in which it is demonstrated that a maximum of 13.8\% gap is guaranteed for problems of max-cut type, and applicable to QCQP problems as the one formulated in this paper.

\bibliographystyle{IEEEtran}

\begin{thebibliography}{10}
\providecommand{\url}[1]{#1}
\csname url@samestyle\endcsname
\providecommand{\newblock}{\relax}
\providecommand{\bibinfo}[2]{#2}
\providecommand{\BIBentrySTDinterwordspacing}{\spaceskip=0pt\relax}
\providecommand{\BIBentryALTinterwordstretchfactor}{4}
\providecommand{\BIBentryALTinterwordspacing}{\spaceskip=\fontdimen2\font plus
\BIBentryALTinterwordstretchfactor\fontdimen3\font minus
  \fontdimen4\font\relax}
\providecommand{\BIBforeignlanguage}[2]{{%
\expandafter\ifx\csname l@#1\endcsname\relax
\typeout{** WARNING: IEEEtran.bst: No hyphenation pattern has been}%
\typeout{** loaded for the language `#1'. Using the pattern for}%
\typeout{** the default language instead.}%
\else
\language=\csname l@#1\endcsname
\fi
#2}}
\providecommand{\BIBdecl}{\relax}
\BIBdecl

\bibitem{Sifuentes2007488}
W.~Sifuentes and A.~Vargas, ``Short-term hydrothermal coordination considering
  an ac network modeling,'' \emph{International Journal of Electrical Power \&
  Energy Systems}, vol.~29, no.~6, pp. 488 -- 496, 2007.

\bibitem{291004}
M.~Piekutowski, T.~Litwinowicz, and R.~Frowd, ``Optimal short-term scheduling
  for a large-scale cascaded hydro system,'' in \emph{Power Industry Computer
  Application Conference, 1993. Conference Proceedings}, May 1993, pp.
  292--298.

\bibitem{1137622}
A.~Conejo, J.~Arroyo, J.~Contreras, and F.~Villamor, ``Self-scheduling of a
  hydro producer in a pool-based electricity market,'' \emph{Power Systems,
  IEEE Transactions on}, vol.~17, no.~4, pp. 1265--1272, Nov 2002.

\bibitem{918302}
T.~Siu, G.~Nash, and Z.~Shawwash, ``A practical hydro, dynamic unit commitment
  and loading model,'' \emph{Power Systems, IEEE Transactions on}, vol.~16,
  no.~2, pp. 301--306, 2001.

\bibitem{IIASA}
M.~P. Nowak, ``A fast descent method for the hydro storage subproblem in power
  generation,'' \emph{Working Paper WP96109 IIASA}, 1996.

\bibitem{589675}
C.~Li, A.~Svoboda, C.-L. Tseng, R.~Johnson, and E.~Hsu, ``Hydro unit commitment
  in hydro-thermal optimization,'' \emph{Power Systems, IEEE Transactions on},
  vol.~12, no.~2, pp. 764--769, 1997.

\bibitem{4562139}
A.~Borghetti, C.~D'Ambrosio, A.~Lodi, and S.~Martello, ``An {MILP} approach for
  short-term hydro scheduling and unit commitment with head-dependent
  reservoir,'' \emph{Power Systems, IEEE Transactions on}, vol.~23, no.~3, pp.
  1115--1124, 2008.

\bibitem{1626389}
E.~Finardi and E.~da~Silva, ``Solving the hydro unit commitment problem via
  dual decomposition and sequential quadratic programming,'' \emph{Power
  Systems, IEEE Transactions on}, vol.~21, no.~2, pp. 835--844, 2006.

\bibitem{1178813}
A.~Borghetti, A.~Frangioni, F.~Lacalandra, and C.~A. Nucci, ``Lagrangian
  heuristics based on disaggregated bundle methods for hydrothermal unit
  commitment,'' \emph{Power Systems, IEEE Transactions on}, vol.~18, no.~1, pp.
  313--323, 2003.

\bibitem{Catalao2010904}
J.~Catal{\~a}o, H.~Pousinho, and V.~Mendes, ``Scheduling of head-dependent
  cascaded reservoirs considering discharge ramping constraints and start/stop
  of units,'' \emph{International Journal of Electrical Power \& Energy
  Systems}, vol.~32, no.~8, pp. 904 -- 910, 2010.

\bibitem{6919349}
M.~Paredes, L.~Martins, and S.~Soares, ``Using semidefinite relaxation to solve
  the day-ahead hydro unit commitment problem,'' \emph{Power Systems, IEEE
  Transactions on}, vol.~30, no.~5, pp. 2695--2705, 2014.

\bibitem{anjos:2012}
M.~F. Anjos and J.~B. Lasserre, ``Introduction to semidefinite, conic and
  polynomial optimization,'' in \emph{Handbook on Semidefinite, Conic, and
  Polynomial Optimization}, ser. International Series in Operations Research \&
  Management Science, M.~F. Anjos and J.~B. Lasserre, Eds.\hskip 1em plus 0.5em
  minus 0.4em\relax Springer, 2012, vol. 166, pp. 1--24.

\bibitem{nesterov:2000}
Y.~Nesterov, H.~Wolkowicz, and Y.~Ye, ``Semidefinite programming relaxations of
  nonconvex quadratic optimization,'' in \emph{Handbook of Semidefinite
  Programming}, ser. International Series in Operations Research \& Management
  Science, H.~Wolkowicz, R.~Saigal, and L.~Vandenberghe, Eds.\hskip 1em plus
  0.5em minus 0.4em\relax Springer, 2000, vol.~27, pp. 361--419.


\bibitem{1997geometry}
M.~Laurent, \emph{Geometry of Cuts and Metrics}, ser. Algorithms and
  Combinatorics.\hskip 1em plus 0.5em minus 0.4em\relax Springer, 1997.

\bibitem{barvinok2002course}
A.~Barvinok, \emph{A Course in Convexity}, ser. Graduate studies in
  mathematics.\hskip 1em plus 0.5em minus 0.4em\relax American Mathematical
  Society, 2002.

\bibitem{7065336}
R.~Madani, M.~Ashraphijuo, and J.~Lavaei, ``Promises of conic relaxation for
  contingency-constrained optimal power flow problem,'' \emph{IEEE Transactions
  on Power Systems}, vol.~31, no.~2, pp. 1297--1307, 2016.

\bibitem{crawford1983}
C.~R. Crawford and Y.~S. Moon, ``Finding a positive definite linear combination
  of two hermitian matrices,'' \emph{Linear Algebra and Its Applications},
  vol.~51, pp. 37 -- 48, 1983.

\bibitem{sdpa}
Makoto Yamashita, Katsuki Fujisawa, Kazuhide Nakata, Maho Nakata, Mituhiro
  Fukuda, Kazuhiro Kobayashi, and Kazushige Goto,
\newblock ``A high-performance software package for semidefinite programs: Sdpa
  7'',
\newblock in {\em Research Report B-460 Dept. of Mathematical and Computing
  Science}, Tokyo~Institute of~Technology, Ed. Tokyo Institute of Technology,
  2010.

\bibitem{Goemans:1995:IAA:227683.227684}
M.~X. Goemans and D.~P. Williamson, ``Improved approximation algorithms for
  maximum cut and satisfiability problems using semidefinite programming,''
  \emph{J. ACM}, vol.~42, no.~6, pp. 1115--1145, Nov. 1995.

\end{thebibliography}

% that's all folks
\end{document}